\documentclass[reqno, 12
pt]{amsart}
\usepackage[a4paper,margin=1in]{geometry}
\usepackage{amsmath,amssymb,amsthm,mathtools}
\usepackage{mathrsfs}
\usepackage{enumitem}
\usepackage{hyperref}
\usepackage{microtype}

\newtheorem{theorem}{Theorem}[section]

\newtheorem{lemma}[theorem]{Lemma}
\newtheorem{definition}[theorem]{Definition}
\newtheorem{corollary}[theorem]{Corollary}
\theoremstyle{remark}
\newtheorem{remark}[theorem]{Remark}
\newtheorem{example}[theorem]{Example}
\newcommand{\Ker}{\text{Ker}}

\title{New conditions for Lipschitz shadowing of nonautonomous dynamics}
\author[D. Dragi\v cevi\'c]{Davor Dragi\v cevi\'c}
\address{Faculty of Mathematics, University of Rijeka, Radmile Matej\v ci\' c 2, 51000 Rijeka, Croatia}
\email{ddragicevic@math.uniri.hr}
\author[M. Onitsuka]{Masakazu Onitsuka}
\address{Department of Applied Mathematics, Faculty of Science, Okayama University of Science, Ridai-cho 1-1, Okayama 700-0005, Japan}
\email{onitsuka@ous.ac.jp}

\date{}

\begin{document}
\maketitle

\begin{abstract}
It is well known that a semilinear ordinary differential equation exhibits Lipschitz shadowing provided that its linear part admits an exponential dichotomy and the nonlinear part is Lipschitz small uniformly in time. The main objective of the present paper is to formulate new conditions for Lipschitz shadowing of semilinear equations that do not require the aforementioned Lipschitz property of the nonlinear term. These new conditions are stated in terms of variational equations associated with the original equation and its pseudosolutions. In the finite-dimensional case, they imply that the associated variational equations exhibit a uniform summable dichotomy. We also obtain corresponding results in the discrete-time case. 

\end{abstract}

\section{Introduction}
Throughout this paper, $\mathbb R$ denotes the set of all real numbers, and $\mathbb R^+ = [0,\infty)$. Let $A\colon \mathbb R^+\to \mathbb R^{n\times n}$ and $f\colon \mathbb R^+\times \mathbb R^n \to \mathbb R^n$ be continuous maps, and consider the associated semilinear differential equation
\begin{equation}\label{eqintro}
    x'(t)=A(t)x(t)+f(t, x(t)), \quad t\in \mathbb R^+.
\end{equation}
We recall that~\eqref{eqintro} is said to be Lipschitz shadowable (or Hyers-Ulam stable) if there exists $L>0$
such that for each $\varepsilon>0$ and a continuously differentiable function $y\colon \mathbb R^+\to \mathbb R^n$ satisfying
\[
\sup_{t\in \mathbb R^+}\|y'(t)-A(t)y(t)-f(t, y(t))\|\le \varepsilon,
\]
there exists a solution $x\colon \mathbb R^+\to \mathbb R^n$ of~\eqref{eqintro} such that
\begin{equation}\label{14:26intro}
\sup_{t\in \mathbb R^+}\|x(t)-y(t)\| \le L\varepsilon.
\end{equation}
This property ensures that in a vicinity of each pseudosolution of~\eqref{eqintro} we can construct its exact solution. This notion matters because each numerical method used to solve~\eqref{eqintro} produces a pseudosolution due to numerical errors. Thus, it is important to know that a true solution exists close to a pseudosolution. We note that this concept is motivated by the notion of shadowing originally introduced and studied for smooth dynamical systems (see~\cite{Palmer, Pil, PilSak}).

It follows from the results in~\cite{BD, BD2} that~\eqref{eqintro} is Lipschitz shadowable provided that the following hold:
\begin{enumerate}
\item[(a)] the linear part of~\eqref{eqintro} given by 
\begin{equation}\label{ldeintro}
x'=A(t)x, \quad t\in \mathbb R^+
\end{equation}
admits an exponential dichotomy (see Remark~\ref{ED});
\item[(b)] for $t\in \mathbb R^+$ and $x, y\in \mathbb R^n$,
\begin{equation}\label{lipintro}
\|f(t, x)-f(t, y)\| \le c\|x-y\|,
\end{equation}
where $c>0$ is sufficiently small.
\end{enumerate}
In addition, the requirement (a) can be slightly weakened by requiring that~\eqref{ldeintro} admits a summable dichotomy (see Definition~\ref{SD}). Furthermore, there are converse results. For example, it turns out that the Lipschitz shadowing of~\eqref{ldeintro} is equivalent to a summable dichotomy (see the proof of~\cite[Theorem 3]{D}). For earlier related results that deal with the Lipschitz shadowing of~\eqref{ldeintro} (or its discrete-time counterpart), we refer to~\cite{BBT, group, BLO, BRST} and references therein. Also relevant to the content of the present paper are the results from~\cite{BDOP}, where the authors studied the so-called conditional shadowing. The main objective was to replace the global Lipschitz condition~\eqref{lipintro} with its local version, showing that this is sufficient to approximate certain (but not necessarily all) pseudosolutions of~\eqref{eqintro} by exact solutions.

The main objective of the present paper is to formulate new conditions for the Lipschitz shadowing of~\eqref{eqintro} (replacing $\mathbb R^n$ with an arbitrary Banach space $X$) which do not require that $f$ be Lipschitz in the second variable. 
These conditions (see Theorem~\ref{T1}) are stated in terms of variational equations related to~\eqref{eqintro} of the form
\begin{equation}\label{varintro}
x'(t)=(A(t)+D_xf(t, q(t)))x(t), \quad t\in \mathbb R^+,
\end{equation}
where $q$ is a bounded perturbation of some pseudosolution of~\eqref{eqintro}. In particular, when $X=\mathbb R^n$, these conditions imply that~\eqref{varintro} admits a summable dichotomy for each $q$. Example~\ref{EX1} illustrates the scope of Theorem~\ref{T1}
by exhibiting a situation in which $f$ is not globally
Lipschitz in its second variable and which falls beyond
the scope of all the previously mentioned results.

We emphasize that all previous results on Lipschitz shadowing
for~\eqref{eqintro} were established using the Banach
contraction principle. In contrast, the proof of Theorem~\ref{T1} employs a fundamentally
different approach, in which $x$ satisfying~\eqref{14:26intro}
is constructed by solving an auxiliary ordinary differential
equation on the Banach space $C_b(\mathbb R^+, X)$ of all
bounded continuous functions from $\mathbb R^+$ to $X$. In Section~\ref{disctime} we establish an analogous result for the case of discrete time.

\section{The case of continuous time}
\subsection{Preliminaries}\label{Prel}
Throughout this paper, $X=(X, \|\cdot \|)$ will be a Banach space and $\mathcal L(X)$ the space of all bounded linear operators on $X$. We will also use the symbol $\| \cdot \|$ to denote the operator norm on $\mathcal L(X)$.  By $C_b(\mathbb R^+, X)$ we will denote the Banach space of all bounded continuous maps $u\colon \mathbb R^+\to X$ equipped with the supremum norm $\|u\|_\infty:=\sup_{t\in \mathbb R^+}\|u(t)\|$.

Take a continuous map $A\colon \mathbb R^+\to \mathcal L(X)$. Furthermore, let $f\colon \mathbb R^+\times X\to X$ be a continuous map such that $x\mapsto f(t,x)$ is continuously Fr\'echet differentiable for each $t\in \mathbb R^+$. Moreover, assume that  $(t,x)\mapsto D_xf(t,x)$ is continuous, where $D_x$ denotes the derivative with respect to $x$. We consider the associated nonautonomous differential equation
\begin{equation}\label{eq:main}
 x'(t)=A(t)x(t)+f(t,x(t)),\qquad t\in \mathbb R^+.
\end{equation}

We recall that~\eqref{eq:main} is \emph{Lipschitz shadowable} (or \emph{Hyers-Ulam stable}) if there exists $L>0$ such that for every $\varepsilon>0$ and a continuously differentiable function $y\colon \mathbb R^+\to X$ such that 
\begin{equation}\label{pseudo}
\sup_{t\in \mathbb R^+}\|y'(t)-A(t)y(t)-f(t, y(t))\|\le \varepsilon,
\end{equation}
there exists a solution $x\colon \mathbb R^+\to X$ of~\eqref{eq:main} satisfying
\begin{equation}\label{dev}
\sup_{t\in \mathbb R^+}\|x(t)-y(t)\| \le L\varepsilon.
\end{equation}

To state the first main result of our paper, we introduce some additional notation. More precisely, for $\varepsilon >0$ let $\mathcal P(\varepsilon)$ consist of all continuously differentiable functions $y\colon \mathbb R^+\to X$ satisfying~\eqref{pseudo}. Furthermore, set $\mathcal P:=\bigcup_{\varepsilon>0}\mathcal P(\varepsilon)$ and 
\[
\mathcal S:=\mathcal P+C_b(\mathbb R^+, X)=\{y+u: \ y\in \mathcal P, \ u\in C_b(\mathbb R^+, X)\}.
\]
\subsection{Main result}
The following is our first main result. It formulates new conditions under which~\eqref{eq:main} is Lipschitz shadowable.
\begin{theorem}\label{T1}
Assume that there exists a family $\{\Gamma_q: \ q\in \mathcal S\}$ of bounded linear operators $\Gamma_q \colon C_b(\mathbb R^+, X)\to C_b(\mathbb R^+, X)$ satisfying the following assumptions:
\begin{enumerate}
\item[(i)] for each $q\in \mathcal S$ and $y\in C_b(\mathbb R^+, X)$, the function $x:=\Gamma_q y$ is differentiable on $\mathbb R^+$ and satisfies
\begin{equation}\label{variational}
x'(t)=A(t)x(t)+D_xf(t, q(t))x(t)+y(t), \quad t\in \mathbb R^+;
\end{equation}
\item[(ii)] 
\begin{equation}\label{uniform}
K:=\sup_{q\in \mathcal S}\|\Gamma_q\|<\infty;
\end{equation}
\item[(iii)] for every $q\in \mathcal S$ and $R>0$, there is a constant $C_{q, R}>0$ such that 
\begin{equation}\label{Lip}
\|\Gamma_{q+z_1}-\Gamma_{q+z_2}\| \le C_{q, R}\|z_1-z_2\|_\infty 
\end{equation}
for every $z_i\in C_b(\mathbb R^+, X)$, $i=1, 2$ such that $\|z_i\|_\infty \le R$ for $i=1, 2$.
\end{enumerate}
Then~\eqref{eq:main} is Lipschitz shadowable.
\end{theorem}

\begin{proof}
Take $y\in \mathcal P(\varepsilon)$ for some $\varepsilon >0$. We define $r\colon \mathbb R^+\to X$ by
\begin{equation}\label{r} 
r(t):=y'(t)- A(t)y(t)-f(t, y(t)), \quad t\in \mathbb R^+.
\end{equation}
Since $y\in \mathcal P(\varepsilon)$, $r\in C_b(\mathbb R^+, X)$ and $\|r\|_\infty \le \varepsilon$. We consider the following autonomous differential equation on $C_b(\mathbb R^+, X)$:
\begin{equation}\label{newode}
\frac{\partial z_\tau}{\partial \tau}=V_y(z_\tau), \quad \tau \in \mathbb R^+,
\end{equation}
where $V_y\colon C_b(\mathbb R^+, X)\to C_b(\mathbb R^+, X)$ is given by
\[
V_y(z):=-\Gamma_{y+z}r, \quad z\in C_b(\mathbb R^+, X).
\]
A solution of~\eqref{newode} is a differentiable function $\mathbb R^+\ni \tau \mapsto z_\tau \in C_b(\mathbb R^+, X)$ which satisfies~\eqref{newode}.

Observe that~\eqref{uniform} implies that $\|V_y(z)\|_\infty \le K\|r\|_\infty \le K\varepsilon$ for each $z\in C_b(\mathbb R^+, X)$. Moreover, by~\eqref{Lip} for $R>0$ there is $C_{y,R}>0$ such that 
\[
\|V_y(z_1)-V_y(z_2)\|_\infty \le C_{y,R} \varepsilon \|z_1-z_2\|_\infty, \quad \text{for $z_i\in C_b(\mathbb R^+, X)$ with $\|z_i\|_\infty \le R$, $i=1, 2$.}
\]
Hence, $V_y$ is locally Lipschitz.

We now need an auxiliary lemma.
\begin{lemma}\label{AL}
The ordinary differential equation~\eqref{newode} has a unique $C^1$ solution $\tau \mapsto z_\tau$ defined globally on $\mathbb R^+$ and satisfying $z_0\equiv 0$. Moreover, 
\begin{equation}\label{additional}
\|z_\tau\|_\infty \le K\varepsilon \tau , \quad \tau \in \mathbb R^+.
\end{equation}
\end{lemma}

\begin{proof}[Proof of the lemma]
Since $V_y$ is locally Lipschitz on $C_b(\mathbb R^+, X)$, the Banach-space version of the Picard--Lindel\"of theorem (see~\cite[Theorem 4.2.4]{AH}) gives a unique maximal solution $\tau\mapsto z_\tau$ of the Cauchy problem consisting of~\eqref{newode} with the initial condition $z_0=0$
on $[0,\tau_*)$ for some $\tau_*\in (0, \infty]$.  From \eqref{uniform}, we get 
\[
\left \|\frac{\partial z_\tau}{\partial \tau} \right \|_\infty \le K\varepsilon,
\]
which implies that 
\begin{equation}\label{611}
\|z_{\tau_1}-z_{\tau_2}\|_\infty=\left \|\int_{\tau_2}^{\tau_1} \frac{\partial z_\tau}{\partial \tau}\, d\tau \right \|_\infty  \le K\varepsilon |\tau_1-\tau_2|, \quad \text{for $\tau_1, \tau_2 \in [0, \tau_*)$.}
\end{equation}
Assume that $\tau_*<\infty$.
By~\eqref{611} and the completeness of $C_b(\mathbb R^+, X)$, $z_*=\lim_{\tau\nearrow \tau_*}z_\tau \in C_b(\mathbb R^+, X)$ exists. By the Picard-Lindel\"of theorem, there exists a continuously differentiable map $\tau \mapsto \bar z_\tau \in C_b(\mathbb R^+, X)$ defined on $[\tau_*, \tau_*+\delta)$ for some $\delta>0$ such that 
\[
\frac{\partial \bar z_\tau}{\partial \tau}=V_y(\bar z_\tau) \quad \text{on $[\tau_*, \tau_*+\delta)$} \quad \text{and} \quad \bar z_{\tau_*}=z_*.
\]
Then
\[
\tilde z_\tau:=\begin{cases}
z_\tau & \tau \in [0, \tau_*); \\
\bar z_\tau & \tau \in [\tau_*, \tau_*+\delta)
\end{cases}
\]
is a solution of~\eqref{newode} on $[0, \tau_*+\delta)$ such that $\tilde z_0=z_0=0$. This contradicts the maximality of $\tau_*$. We conclude that $\tau_*=\infty$.
\end{proof}
Let $\tau \mapsto z_\tau$ be given by Lemma~\ref{AL}. For $\tau \in \mathbb R^+$, set
\[
u_\tau:=y+z_\tau \quad \text{and} \quad w_\tau:=\Gamma_{u_\tau}r.
\]
Then $\frac{\partial u_\tau}{\partial \tau}=V_{y}(z_\tau)=-w_\tau$. Moreover, by~\eqref{variational} we have 
\begin{equation}\label{653}
w_\tau'(t)-A(t)w_\tau(t)-D_xf(t, u_\tau(t))w_\tau(t)=r(t), \quad \text{for every $t, \tau \in \mathbb R^+$.} 
\end{equation}

\begin{lemma}
For every $\tau \in \mathbb R^+$,
\begin{equation}\label{lemeq}
u_\tau'(t)-A(t)u_\tau(t)-f(t, u_\tau(t))=(1-\tau)r(t), \quad t\in \mathbb R^+.
\end{equation}
\end{lemma}

\begin{proof}[Proof of the lemma]
Since $u_0=y$, integrating $\frac{\partial u_\tau}{\partial \tau}=-w_\tau$ gives that for $\tau \in \mathbb R^+$,
\[
u_\tau(t)=y(t)-\int_0^\tau w_s(t)\, ds, \quad t\in \mathbb R^+.
\]
We now claim that for $\tau \in \mathbb R^+$,
\begin{equation}\label{821}
u_\tau'(t)=y'(t)-\int_0^\tau w_s'(t)\, ds, \quad t\in \mathbb R^+.
\end{equation}
To this end, it is sufficient to show that the maps $(s, t)\mapsto w_s(t)$ and $(s, t)\mapsto w_s'(t)$ are continuous.  Notice that  it follows from~\eqref{Lip} and~\eqref{611} that for arbitrary $M>0$ and a compact interval $J\subset \mathbb R^+$, there exists $C>0$ such that 
\[
\begin{split}
\|w_{s}(t)-w_{s'}(t')\| &\le \|w_{s}(t)-w_{s'}(t)\|+\|w_{s'}(t)-w_{s'} (t')\| \\
&\le \|w_{s}-w_{s'}\|_\infty+\|w_{s'}(t)-w_{s'}(t')\| \\
&\le \|\Gamma_{y+z_{s}}-\Gamma_{y+z_{s'}}\| \cdot \|r\|_\infty+\|w_{s'}(t)-w_{s'}(t')\|  \\
&\le C|s-s'|+\|w_{s'}(t)-w_{s'}(t')\|
\end{split}
\]
for all $(s, t), (s', t')\in [0, M]\times J$. Since $t\mapsto w_s(t)$ is differentiable and thus continuous, we have obtained the continuity of the map $(s, t)\mapsto w_s(t)$. Similarly, using~\eqref{611} and the fact that
$u_s=y+z_s$, we conclude that the map
$(s,t)\mapsto u_s(t)$ is continuous.
Consequently, the continuity of the maps
$A$, $D_xf$ and $r$, together with (16),
yields the continuity of the map
$(s,t)\mapsto w_s'(t)$.
This establishes~\eqref{821}.
By~\eqref{653} and~\eqref{821}, we have
\begin{equation*}
\begin{split}
u_\tau'(t) &=y'(t)-A(t)\int_0^\tau w_s(t)\, ds-\int_0^\tau D_xf(t, u_s(t))w_s(t)\, ds-\tau r(t) \\
&=y'(t)+A(t)(u_\tau(t)-y(t))-\int_0^\tau D_xf(t, u_s(t))w_s(t)\, ds-\tau r(t),
\end{split}
\end{equation*}
and therefore for $\tau \in \mathbb R^+$,
\begin{equation}\label{1006}
u_\tau'(t)-A(t)u_\tau(t)=y'(t)-A(t)y(t)-\int_0^\tau D_xf(t, u_s(t))w_s(t)\, ds-\tau r(t), \quad t\in \mathbb R^+.
\end{equation}
Next, the chain-rule gives
\[
\partial_s f(t, u_s(t))=-D_x f(t, u_s(t))w_s(t), 
\]
and by integrating, we obtain 
\[
-\int_0^\tau D_xf(t, u_s(t))w_s(t)\, ds=f(t, u_\tau(t))-f(t, y(t)).
\]
Inserting this into~\eqref{1006} and using the definition of $r$ (see~\eqref{r}), we obtain~\eqref{lemeq}.
\end{proof}
Set $x:=u_1$. It follows from~\eqref{lemeq} for $\tau=1$ that $x$ is a solution of~\eqref{eq:main}. Moreover, \eqref{additional} gives that 
\[
\|x-y\|_\infty=\|z_1\|_\infty \le K\varepsilon,
\]
yielding that~\eqref{dev} holds with $L:=K>0$. The proof of the theorem is complete.
\end{proof}

\subsection{Relation to the previous results and an explicit example}
We consider the linear part of~\eqref{eq:main} given by
\begin{equation}\label{ldemain}
x'(t)=A(t)x(t), \quad t\in \mathbb R^+.
\end{equation}
By $\{T(t, s): \ t, s\in \mathbb R^+\}$ we will denote the evolution family associated to~\eqref{ldemain}.
We recall the notion of a summable dichotomy.
\begin{definition}\label{SD}
We say that~\eqref{ldemain} admits a \emph{summable dichotomy} if there exist a projection $P\in \mathcal L(X)$ and a constant $K>0$ such that
\begin{equation}\label{sumdic}
\int_0^t \|T(t, 0)PT(0, s)\|\, ds+\int_t^\infty\|T(t, 0)(I-P)T(0, s)\|\, ds\le K, \quad t\in \mathbb R^+,
\end{equation}
where $I$ denotes the identity operator on $X$.
\end{definition}

\begin{remark}\label{ED}
We recall that~\eqref{ldemain} admits an exponential dichotomy if there exist a projection $P\in \mathcal L(X)$ and constants $D, \lambda>0$ such that
\begin{enumerate}
\item for $t\ge s\ge 0$,
\begin{equation}\label{ED1}
\|T(t, 0)PT(0, s)\| \le De^{-\lambda (t-s)};
\end{equation}
\item for $0\le t\le s$,
\begin{equation}\label{ED2}
\|T(t, 0)(I-P)T(0, s)\|\le De^{-\lambda (s-t)}.
\end{equation}
\end{enumerate}
Clearly, the notion of an exponential dichotomy is a particular case of the notion of a summable dichotomy, as~\eqref{ED1} and~\eqref{ED2} imply that~\eqref{sumdic} holds with $K=\frac{2D}{\lambda}$.
\end{remark}

We have the following consequence of Theorem~\ref{T1}.
\begin{corollary}\label{Cor1}
Suppose that~\eqref{ldemain} admits a summable dichotomy, and let $P\in \mathcal L(X)$ and $K>0$ be such that~\eqref{sumdic} holds. Furthermore, suppose that there is $c>0$ satisfying
\begin{equation}\label{c}
c<\frac 1 K
\end{equation}
such that 
\begin{equation}\label{7:23}
\|D_xf(t, x)\| \le c, \quad \text{for $t\in \mathbb R^+$ and $x\in X$.}
\end{equation}
Finally, suppose that there is $L>0$ such that 
\begin{equation}\label{L}
\|D_xf(t, x)-D_xf(t, y)\| \le L\|x-y\|, \quad \text{for $t\in \mathbb R^+$ and $x, y\in X$.}
\end{equation}
 Then~\eqref{eq:main} is Lipschitz shadowable.
\end{corollary}

\begin{proof}
Let $q\colon \mathbb R^+\to X$ be a continuous map and take an arbitrary $y\in C_b(\mathbb R^+, X)$. We consider a map $\mathcal T_q\colon C_b(\mathbb R^+, X)\to C_b(\mathbb R^+, X)$ given by 
\[
(\mathcal T_q x)(t)=\int_0^\infty \mathcal G(t, s)(D_x f(s, q(s))x(s)+y(s))\, ds, \quad t\in \mathbb R^+, \ x\in C_b(\mathbb R^+, X)
\]
where
\[
\mathcal G(t, s)=\begin{cases}
T(t, 0)PT(0, s) & t\ge s; \\
-T(t, 0)(I-P)T(0, s) &t <s.
\end{cases}
\]
By~\eqref{sumdic} and~\eqref{7:23},
\begin{equation}\label{estimate}
\|(\mathcal T_q x)(t)\|\le K(c\|x\|_\infty+\|y\|_\infty)
\end{equation}
for $t\in \mathbb R^+$ and $x\in C_b(\mathbb R^+, X)$. Hence, $\mathcal T_q x\in C_b(\mathbb R^+, X)$ and thus $\mathcal T_q$ is well-defined. Similarly, 
\[
\|\mathcal T_q x_1-\mathcal T_q x_2\|_\infty \le Kc\|x_1-x_2\|_\infty, \quad \text{for $x_1, x_2\in C_b(\mathbb R^+, X)$.}
\]
From~\eqref{c} we conclude that $\mathcal T_q$ is a contraction on $C_b(\mathbb R^+, X)$. Set $\Gamma_q y:=x$, where $x\in C_b(\mathbb R^+, X)$ is the unique fixed point of $\mathcal T_q$. 

To verify~\eqref{variational}, set $\omega(s):=D_x f(s, q(s))x(s)+y(s)$ for $s\in \mathbb R^+$. Since $D_x f$, $q$, $x$, and $y$ are continuous, $\omega$ is continuous, and by~\eqref{7:23} we have $\|\omega(s)\|\le c\|x\|_\infty+\|y\|_\infty$ for all $s\in \mathbb R^+$. In particular, $\omega\in C_b(\mathbb R^+, X)$, and~\eqref{sumdic} with $t=0$ ensures that
\[
\int_0^\infty \|(I-P)T(0, s)\omega(s)\|\, ds \le K\|\omega\|_\infty < \infty.
\]
Since $x=\mathcal T_q x$ and $T(t, 0)\in \mathcal L(X)$ is independent of $s$, we can rewrite $x(t)$ as
\[
x(t)=T(t, 0)\int_0^t PT(0, s)\omega(s)\, ds - T(t, 0)\int_t^\infty (I-P)T(0, s)\omega(s)\, ds, \quad t\in \mathbb R^+.
\]
Since $t\mapsto T(t, 0)$ is continuously differentiable with $\frac{d}{dt}T(t, 0)=A(t)T(t, 0)$ and the integrands are continuous, differentiating with respect to $t$ yields
\[
\begin{split}
x'(t) &= A(t)T(t, 0)\int_0^t PT(0, s)\omega(s)\, ds + T(t, 0)PT(0, t)\omega(t) \\
&\quad - A(t)T(t, 0)\int_t^\infty (I-P)T(0, s)\omega(s)\, ds + T(t, 0)(I-P)T(0, t)\omega(t) \\
&= A(t)x(t) + T(t, 0)(P + I - P)T(0, t)\omega(t) = A(t)x(t) + D_x f(t, q(t))x(t) + y(t)
\end{split}
\]
for every $t\in \mathbb R^+$, where we used the identity $T(t, 0)T(0, t)=I$. Thus, \eqref{variational} holds.

By~\eqref{estimate} we have that 
\[
\|\Gamma_q y\|_\infty = \|x\|_\infty \le \frac{K}{1-Kc}\|y\|_\infty, 
\]
which implies that~\eqref{uniform} holds.

We still need to establish~\eqref{Lip}. To this end, we note that for every continuous map $q\colon \mathbb R^+\to X$ and $z_i\in C_b(\mathbb R^+, X)$, $i=1, 2$, we have (using~\eqref{L}) that 
\[
\|\mathcal T_{q+z_1}x-\mathcal T_{q+z_2}x\|_\infty \le KL\|z_1-z_2\|_\infty \cdot \|x\|_\infty, \quad x\in C_b(\mathbb R^+, X).
\]
Hence, letting $x_i:=\Gamma_{q+z_i}y$, $i=1, 2$, we have
\[
\begin{split}
\|x_1-x_2\|_\infty &=\|\mathcal T_{q+z_1}x_1-\mathcal T_{q+z_2}x_2\|_\infty  \\
&\le \|\mathcal T_{q+z_1}x_1-\mathcal T_{q+z_2}x_1\|_\infty+\|\mathcal T_{q+z_2}x_1-\mathcal T_{q+z_2}x_2\|_\infty \\
&\le KL\|z_1-z_2\|_\infty \cdot \|x_1\|_\infty+Kc\|x_1-x_2\|_\infty \\
&\le KL\cdot \frac{K}{1-Kc}\|z_1-z_2\|_\infty \cdot \|y\|_\infty+Kc\|x_1-x_2\|_\infty.
\end{split}
\]
We conclude that 
\[
\|x_1-x_2\|_\infty \le \frac{K^2L}{(1-Kc)^2}\|z_1-z_2\|_\infty \cdot \|y\|_\infty,
\]
which yields~\eqref{Lip}.
\end{proof}

\begin{remark}
Corollary~\ref{Cor1} is not a new result as it follows from the results in~\cite{BDS} and a simple observation that~\eqref{7:23} implies that
\[
\|f(t, x)-f(t, y)\| \le c\|x-y\|, \quad t\in \mathbb R^+, \ x, y\in X.
\]
\end{remark}

The following example shows that Theorem~\ref{T1} applies beyond Corollary~\ref{Cor1}.
\begin{example}\label{EX1}
Let $X=\mathbb R^2$ with the norm $\|(x, y)\|=\max \{|x|, |y|\}$, $(x, y)\in X$. We consider~\eqref{eq:main} with 
\[
A(t)=\begin{pmatrix}
-1 & 0\\
0 & 1
\end{pmatrix} \quad \text{and} \quad f(t, x, y)=\left (-\frac 12 x|x|+p(t), \frac 1 2 y|y|\right ),
\]
for $t\in \mathbb R^+$ and $(x, y)\in \mathbb R^2$, where $p(t)=1+t+ \frac 1 2t^2$. It is straightforward to verify that $f$ has the regularity properties imposed in Section~\ref{Prel}. In addition, note that 
\begin{equation}\label{derivative}
D_{(x, y)}f(t, x, y)=\begin{pmatrix}
-|x| & 0\\
0 & |y|
\end{pmatrix}, \quad t\in \mathbb R^+, \ (x, y)\in \mathbb R^2.
\end{equation}
In particular, $f$ is not globally Lipschitz in the variable $(x, y)$. Hence, Corollary~\ref{Cor1} is not applicable.
We observe that 
\begin{equation}\label{4:38}
\|D_{(x, y)}f(t, x_1, y_1)-D_{(x, y)}f(t, x_2, y_2)\|\le \|(x_1, y_1)-(x_2, y_2)\|, \quad (x_i, y_i)\in \mathbb R^2, \ i=1, 2.
\end{equation}

Take an arbitrary  continuous map $q=(q_1, q_2)\colon \mathbb R^+\to \mathbb R^2$ and define $\Gamma_q\colon C_b(\mathbb R^+, \mathbb R^2)\to C_b(\mathbb R^+, \mathbb R^2)$ by 
\[
(\Gamma_q h)(t):=\left (\int_0^t e^{-\int_s^t (1+|q_1(\xi)|)\, d\xi}h_1(s)\, ds, -\int_t^\infty e^{-\int_t^s(1+|q_2(\xi)|)\, d\xi}h_2(s)\, ds\right ),
\]
for $t\in \mathbb R^+$ and $h=(h_1, h_2)\in C_b(\mathbb R^+, \mathbb R^2)$. Note that $\Gamma_q$ is a linear operator on $C_b(\mathbb R^+, \mathbb R^2)$. Moreover, it is bounded as 
\begin{equation}\label{303}
\|\Gamma_q h\|_\infty \le \|h\|_\infty, \quad h\in C_b(\mathbb R^+, \mathbb R^2).
\end{equation}
A simple computation shows that for each $y\in C_b(\mathbb R^+, \mathbb R^2)$, $x=\Gamma_q y$ is a solution of~\eqref{variational}.  In addition, \eqref{303} implies that~\eqref{uniform} holds.

Take now $q\in \mathcal S$ and $z_i\in C_b(\mathbb R^+, \mathbb R^2)$, $i=1, 2$. Moreover, take an arbitrary $h\in C_b(\mathbb R^+, \mathbb R^2)$ and set $x_i:=\Gamma_{q+z_i}h\in C_b(\mathbb R^+, \mathbb R^2)$, $i=1, 2$, $x:=x_1-x_2\in C_b(\mathbb R^+, \mathbb R^2)$. Then
\begin{equation}\label{4:41}
\begin{split}
x'(t) &=A(t)x(t)+D_{(x, y)}f(t, q(t)+z_1(t))x(t) \\
&\phantom{=}+\left [D_{(x, y)}f(t, q(t)+z_1(t))-D_{(x, y)}f(t, q(t)+z_2(t))\right ]x_2(t), \quad t\in \mathbb R^+.
\end{split}
\end{equation}
By~\eqref{4:38}, 
\[
\left [D_{(x, y)}f(\cdot, q(\cdot)+z_1(\cdot))-D_{(x, y)}f(\cdot, q(\cdot)+z_2(\cdot))\right ]x_2(\cdot)\in C_b(\mathbb R^+, \mathbb R^2).
\]
Since $x\in C_b(\mathbb R^+, \mathbb R^2)$ and $x(0)\in \{0\}\times \mathbb R$, \eqref{4:41} gives that 
\begin{equation}\label{2135}
x=\Gamma_{q+z_1}\left (\left [D_{(x, y)}f(\cdot, q(\cdot)+z_1(\cdot))-D_{(x, y)}f(\cdot, q(\cdot)+z_2(\cdot))\right ]x_2(\cdot)\right ).
\end{equation}
Indeed,  the equation
\[
w'(t)=A(t)w(t)+D_{(x, y)}f(t, q(t)+z_1(t))w(t), \quad t\in \mathbb R^+ 
\]
admits an exponential dichotomy with respect to the projection
\[
P=\begin{pmatrix}
1 & 0\\
0 & 0
\end{pmatrix}.
\]
To verify this assertion, write $q+z_1=(v_1,v_2)$. By~\eqref{derivative},
the coefficient matrix of the preceding equation is
\[
B(t)=\begin{pmatrix}-1-|v_1(t)|&0\\0&1+|v_2(t)|\end{pmatrix}.
\]
Consequently, its evolution family is given by
\[
U(t,s)=\begin{pmatrix}
\displaystyle\exp\left(-\int_s^t(1+|v_1(\xi)|)\,d\xi\right)&0\\
0&\displaystyle\exp\left(\int_s^t(1+|v_2(\xi)|)\,d\xi\right)
\end{pmatrix}, \quad t, s\in \mathbb R^+.
\]
For $t\ge s\ge0$,
\[
\|U(t,0)PU(0, s)\|=\exp\left(-\int_s^t(1+|v_1(\xi)|)\,d\xi\right)
\le e^{-(t-s)},
\]
whereas for $0\le t\le s$,
\[
\|U(t,0)(I-P)U(0, s)\|=\exp\left(-\int_t^s(1+|v_2(\xi)|)\,d\xi\right)
\le e^{-(s-t)}.
\]
Thus all the conditions in the definition of exponential dichotomy are
satisfied with $D=\lambda=1$, independently of the choice of $q$ and $z_1$.

Consequently, for each $\rho \in C_b(\mathbb R^+, \mathbb R^2)$ there exists a unique $w\in C_b(\mathbb R^+, \mathbb R^2)$ with $w(0)\in \{0\}\times \mathbb R=\Ker P$ and 
\[
w'(t)=A(t)w(t)+D_{(x, y)}f(t, q(t)+z_1(t))w(t)+\rho(t), \quad t\in \mathbb R^+.
\]
This implies~\eqref{2135}.
Consequently,
\[
\begin{split}
\|\Gamma_{q+z_1}h-\Gamma_{q+z_2}h\|_\infty &=\|x\|_\infty \\
&=\left \|\Gamma_{q+z_1}\left (\left [D_{(x, y)}f(\cdot, q(\cdot)+z_1(\cdot))-D_{(x, y)}f(\cdot, q(\cdot)+z_2(\cdot))\right ]x_2(\cdot)\right )\right \|_\infty \\
&\le \left \| \left [D_{(x, y)}f(\cdot, q(\cdot)+z_1(\cdot))-D_{(x, y)}f(\cdot, q(\cdot)+z_2(\cdot))\right ]x_2(\cdot)\right \|_\infty \\
&\le \|z_1-z_2\|_\infty \cdot \|x_2\|_\infty \\
&=\|z_1-z_2\|_\infty \cdot \|\Gamma_{q+z_2}h\|_\infty \\
&\le \|z_1-z_2\|_\infty \cdot \|h\|_\infty,
\end{split}
\]
where we used~\eqref{4:38} and~\eqref{303}. We conclude that~\eqref{Lip} holds with $C_{q, R}=1$. Therefore, Theorem~\ref{T1} is applicable.

We now argue that the previous results on conditional Lipschitz shadowing in~\cite{BDOP} do not apply to this example. To this end, 
take $\delta \in [0, 1]$ and define $y_\delta\colon \mathbb R^+\to \mathbb R^2$ by
\[
y_\delta(t)=\left (t+\frac{\delta}{1+t^2}, 0\right ),  \quad t\in \mathbb R^+.
\]
Note that $y_0$ is a solution of~\eqref{eq:main}. Let $\eta_\delta(t)=\frac{\delta}{1+t^2}$ so that $y_\delta=y_0+(\eta_\delta,0)$. Write $\psi(s)=\frac12s^2$, so $\psi'(s)=s$. Note that 
\[
r_\delta(t):=y_\delta'(t)-A(t)y_\delta(t)-f(t, y_\delta(t))=\left(\eta_\delta'(t)+\eta_\delta(t)+\psi(t+\eta_\delta(t))-\psi(t),\,0\right).
\]
Since $\psi(t+\eta)-\psi(t)=\int_0^\eta (t+s)\,ds$ for $\eta\in \mathbb R^+$,
\[
0\le \psi(t+\eta)-\psi(t)= t\eta+\tfrac12\eta^2.
\]
In addition,
\[
\|\eta_\delta\|_\infty=\delta,\quad
\sup_{t\in \mathbb R^+}(t\eta_\delta(t))=\frac\delta2,\quad
\|\eta_\delta'\|_\infty
=\delta\sup_{t\in \mathbb R^+}\frac{2t}{(1+t^2)^2}\le\delta.
\]
Combining these bounds, we get
\[
\|r_\delta\|_\infty\le 
\frac 5 2 \delta+\frac12\delta^2.
\]
Note that $r_\delta(0)=(\delta+\frac 1 2 \delta^2, 0)$, which means that $y_\delta$ is a proper pseudosolution of~\eqref{eq:main} for $\delta\neq 0$.

Take $\delta \in (0, 1]$. We now claim that there does not exist an open set $H\subset \mathbb R^2$ with the following properties:
\begin{enumerate}
\item[(i)] $y_\delta(t)\in H$ for $t\in \mathbb R^+$;
\item[(ii)] $f(t, \cdot)$ is uniformly (in $t$) Lipschitz on $H$. 
\end{enumerate}
Indeed, otherwise we would have that 
\[
\sup_{t\in \mathbb R^+}\|D_{(x, y)}f(t, y_\delta(t))\|<\infty, 
\]
which is impossible due to~\eqref{derivative} and the fact that $\sup_{t\in \mathbb R^+}\|y_\delta(t)\|=\infty$. Hence, \cite[Theorem~2]{BDOP} is not applicable. 

Finally, denoting
\[
g(t, x, y)=A(t)(x, y)+f(t, x, y), 
\]
for the logarithmic norm (see~\cite[Section 3]{BDOP}) $\mu_\infty$, we have 
\[
\mu_\infty(D_{(x, y)}g(t, x, y))=1+|y|\ge 1, 
\]
which means that~\cite[Theorem 5]{BDOP} is not applicable.

\end{example}

\begin{remark}
Assume that $X$ is finite-dimensional. The assumptions (i) and (ii) from the statement of Theorem~\ref{T1} together with~\cite[Theorem 1, p.131]{Coppel} (and its proof) imply that for each $q\in \mathcal S$, the equation
\[
x'(t)=(A(t)+D_xf(t, q(t)))x(t), \quad t\in \mathbb R^+
\]
admits a summable dichotomy. Moreover, the associated constant $K>0$ in~\eqref{sumdic} can be chosen uniformly over all $q$.
\end{remark}

\section{The case of discrete time}\label{disctime}
We now present a discrete-time version of Theorem~\ref{T1}. Throughout this section, let $\mathbb Z^+:=\{0, 1, 2, \dots\}$.

Let $(A_n)_{n \in \mathbb Z^+}$ be a sequence in
$\mathcal L(X)$, and let $f_n:X\to X$ be continuously Fr\'echet
differentiable for every $n \in \mathbb Z^+$. We do not assume that any $A_n$ is
invertible. Consider the difference equation
\begin{equation}\label{eq:discrete}
 x_{n+1}=A_nx_n+f_n(x_n),\qquad n \in \mathbb Z^+.
\end{equation}
Write $\ell^\infty=\ell^\infty(\mathbb Z^+,X)$, with its usual
supremum norm. We say that~\eqref{eq:discrete} is \emph{Lipschitz shadowable} if there exists $L>0$ such that for every $\varepsilon>0$ and a sequence $(y_n)_{n \in \mathbb Z^+}\subset X$ such that 
\begin{equation}\label{pseudo-d}
\sup_{n \in \mathbb Z^+}\|y_{n+1}-A_n y_n-f_n(y_n)\|\le \varepsilon,
\end{equation}
there exists a solution $(x_n)_{n \in \mathbb Z^+}$ of~\eqref{eq:discrete} satisfying
\begin{equation}\label{dev-d}
\sup_{n \in \mathbb Z^+}\|x_n-y_n\| \le L\varepsilon.
\end{equation}

For $\varepsilon >0$ let $\mathcal P(\varepsilon)$ consist of all sequences $(y_n)_{n \in \mathbb Z^+}\subset X$ satisfying~\eqref{pseudo-d}. Furthermore, set $\mathcal P:=\bigcup_{\varepsilon>0}\mathcal P(\varepsilon)$ and 
\[
\mathcal S:=\mathcal P+\ell^\infty.
\]

The following is the second main result of our paper.
\begin{theorem}\label{thm:discrete}
There exists a family $\{\Gamma_q: \ q\in \mathcal S\}$ of bounded linear operators $\Gamma_q:\ell^\infty\to\ell^\infty$ satisfying the following assumptions:
\begin{enumerate}
\item[(i)] for each $h\in\ell^\infty$, $v:=\Gamma_qh$ satisfies
\begin{equation}\label{eq:discrete-var}
 v_{n+1}=\bigl(A_n+Df_n(q_n)\bigr)v_n+h_n,
 \qquad n \in \mathbb Z^+;
\end{equation}
\item[(ii)] $K:=\sup_{q\in\mathcal S}\|\Gamma_q\|<\infty$;
\item[(iii)] for every $q\in\mathcal S$ and $R>0$, there is a constant
$C_{q,R}>0$ such that
\[
 \|\Gamma_{q+z_1}-\Gamma_{q+z_2}\|
 \le C_{q,R}\|z_1-z_2\|_\infty
\]
for every $z_i\in\ell^\infty$, $i=1, 2$ such that $\|z_i\|_\infty \le R$ for $i=1, 2$.
\end{enumerate}
Then \eqref{eq:discrete} is Lipschitz shadowable.
\end{theorem}

\begin{proof}
Take $y=(y_n)_{n \in \mathbb Z^+}\in  \mathcal P(\varepsilon)$ for some $\varepsilon >0$, and put
$r_n=y_{n+1}-A_ny_n-f_n(y_n)$, so that
$r=(r_n)_{n \in \mathbb Z^+}\in\ell^\infty$ and $\|r\|_\infty\le\varepsilon$.

By the same reasoning as in the proof of Theorem~\ref{T1}, there exists a unique $C^1$ map $\mathbb R^+\ni \tau \mapsto z_\tau \in \ell^\infty$ satisfying 
\begin{equation}\label{eq:discrete-lift}
 \frac{dz_\tau}{d\tau}=-\Gamma_{y+z_\tau}r,
 \qquad z_0=0.
\end{equation}
Moreover, integrating this and using (ii), we get
\[\|z_\tau\|_\infty\le K\tau \varepsilon, \quad \tau \in \mathbb R^+.\]

Set $u_\tau=y+z_\tau$ and $w_\tau=\Gamma_{u_\tau}r$.
Then $\frac{d u_\tau}{d\tau} =-w_\tau$ and due to~\eqref{eq:discrete-var},
\[
 (w_\tau)_{n+1}
 =\bigl(A_n+Df_n((u_\tau)_n)\bigr)(w_\tau)_n+r_n, \quad n \in \mathbb Z^+.
\]
For each fixed $n \in \mathbb Z^+$, the map
\[
 F_n(\tau)=(u_\tau)_{n+1}-A_n(u_\tau)_n
                   -f_n((u_\tau)_n)
\]
is continuously differentiable, and the chain rule yields
\[
 \begin{split}
 F_n'(\tau)
 &=-(w_\tau)_{n+1}+A_n(w_\tau)_n
       +Df_n((u_\tau)_n)(w_\tau)_n\\
 &=-r_n.
 \end{split}
\]
Since $F_n(0)=r_n$, we obtain $F_n(\tau)=(1-\tau)r_n$.
Thus $x=u_1$ is an exact solution of \eqref{eq:discrete}, and
\[
 \|x-y\|_\infty=\|z_1\|_\infty
 \le K\varepsilon.
\]
The proof is complete.
\end{proof}

\begin{remark}
We note that the most general result on Lipschitz shadowing
for~\eqref{eq:discrete} under the assumption that each $f_n$
is globally Lipschitz is given in~\cite[Theorem 3.4]{DP},
although the authors of~\cite{DP} considered~\eqref{eq:discrete}
on $\mathbb Z$. The following example demonstrates the
applicability of Theorem~\ref{thm:discrete} beyond the
globally Lipschitz setting.
\end{remark}

\begin{example}
Let $X=\mathbb R^2$ with the norm $\|(x, y)\|=\max \{|x|, |y|\}$, $(x, y)\in X$. We consider~\eqref{eq:discrete} with 
\[
A_n=0 \quad \text{and} \quad f_n( x, y)=\left (\frac 1 2 x, n+1+(n+2)(y-n)+(y-n)^3 \right )
\]
for $n \in \mathbb Z^+$ and $(x, y)\in \mathbb R^2$. Note that $f_n$ is not Lipschitz.

Take a sequence $q=(q_n)_{n \in \mathbb Z^+}\subset \mathbb R^2$, and write $q_n=(q_n^1, q_n^2)$, $n \in \mathbb Z^+$. We define $\Gamma_q \colon \ell^\infty \to \ell^\infty$ by 
\[
(\Gamma_q h)_n =\left (\sum_{j=0}^{n-1}2^{-(n-1-j)}h_j^1, -\sum_{j=n}^\infty\frac{h_j^2}{b_n(q)\cdots b_j(q)}\right ), \quad n \in \mathbb Z^+, \ h=(h_n)_{n \in \mathbb Z^+}\in \ell^\infty, \ h_n=(h_n^1, h_n^2),
\]
where 
\[
 b_n(q)=n+2+3(q_n^2-n)^2, \quad n \in \mathbb Z^+.
\]
Clearly, $\Gamma_q$ is linear. Moreover, since $b_n(q)\ge 2$ for each $n \in \mathbb Z^+$, it is easy to verify that $\|\Gamma_q\| \le 2$. A straightforward computation shows that for each $h\in\ell^\infty$, $v=\Gamma_qh$ satisfies~\eqref{eq:discrete-var}. Thus, the assumptions $(i)$ and $(ii)$ in Theorem~\ref{thm:discrete} are satisfied. 

Take now two sequences $q=(q_n)_{n \in \mathbb Z^+}\subset \mathbb R^2$ and $\tilde q=(\tilde q_n)_{n \in \mathbb Z^+}\subset \mathbb R^2$ such that $q-\tilde q\in \ell^\infty$. Moreover, choose an arbitrary $h=(h_n)_{n \in \mathbb Z^+}\in \ell^\infty$. Note that $(\Gamma_qh)_n^1=(\Gamma_{\tilde q}h)_n^1$ for each $n \in \mathbb Z^+$. We now want to compare $(\Gamma_qh)_n^2$ and $(\Gamma_{\tilde q}h)_n^2$. To this end, set $\phi_n(s)=\log(n+2+3(s-n)^2)$. Then
\[
 |\phi_n'(s)|=\frac{6|s-n|}{n+2+3(s-n)^2}
 \le\sqrt{\frac3{n+2}}\le\sqrt{\frac32}=:C.
\]
For $j\ge n$ write 
\[
P_{n,j}(q)=\prod_{i=n}^j b_i(q)^{-1}=e^{-\sum_{i=n}^j \phi_i(q_i^2)} \ \quad \text{and} \quad \ P_{n,j}(\tilde q)=\prod_{i=n}^j b_i(\tilde q)^{-1}=e^{-\sum_{i=n}^j \phi_i(\tilde q_i^2)}.
\]
Both exponents $-\sum_{i=n}^j\phi_i(q_i^2)$ and $-\sum_{i=n}^j\phi_i(\tilde q_i^2)$ are at most $-(j-n+1)\log2$. The mean value theorem  yields
\[
 |P_{n,j}(q)-P_{n,j}(\widetilde q)|
 \le 2^{-(j-n+1)}(j-n+1)C\|q-\widetilde q\|_\infty.
\]
Summing and using $\sum_{k\ge1}k2^{-k}=2$ verifies the assumption $(iii)$ in the statement of Theorem~\ref{thm:discrete}. We conclude that Theorem~\ref{thm:discrete} is applicable.

For $\delta \in (0, 1]$, let $z^\delta=(z_n^\delta)_{n \in \mathbb Z^+}\subset \mathbb R^2$ be given by 
\[
z_n^\delta=\left (0, n+\frac{\delta}{(n+1)^2}\right), \quad n \in \mathbb Z^+.
\]
It is easy to verify that $z^\delta$ satisfies~\eqref{pseudo-d} with $\varepsilon=5\delta$ and that $z^\delta$ is not a solution of~\eqref{eq:discrete}. Furthermore, since
\[
Df_n(z_n^\delta)=\begin{pmatrix}
\frac 1 2 & 0\\
0 & n+2+\frac{3\delta^2}{(n+1)^4}
\end{pmatrix}, \quad n \in \mathbb Z^+, 
\]
and thus 
\[
\sup_{n \in \mathbb Z^+}\|Df_n(z_n^\delta)\|=\infty.
\]
Consequently, as in Example~\ref{EX1}, the discrete-time counterparts to the results in~\cite {BDOP} do not apply.
\end{example}

\section{Conclusions}
In this paper, we have developed a general functional-analytic framework for establishing Lipschitz shadowing for both continuous-time and discrete-time semilinear nonautonomous systems. By formulating the shadowing problem in terms of a family of linear solution operators $\Gamma_q$ associated with the variational equations along pseudosolutions, we proved that uniform boundedness and local Lipschitz dependence of $\Gamma_q$ with respect to $q$ are sufficient to guarantee Lipschitz shadowability (Theorems~\ref{T1} and~\ref{thm:discrete}).

A key advantage of our approach is that it does not require the linear part to admit a classical exponential dichotomy, nor does it require the nonlinear term to be globally Lipschitz or the linear operators $A_n$ in the discrete-time setting to be invertible. In particular, we showed that our results apply to systems whose linear part admits only a summable dichotomy (Corollary~\ref{Cor1}), as well as to explicit nonlinear systems with unbounded trajectories where the nonlinearity fails to be globally Lipschitz and existing criteria such as~\cite[Theorem~2]{BDOP} are not applicable.

\section*{AI disclosure statement}
 ChatGPT (GPT-6.2) was used to suggest the novel approach
to Lipschitz shadowing developed in this work. All proofs
were established independently by the authors.

	\section*{Funding declaration}  
This paper has been funded by the European Union – NextGenerationEU-StatDin\-MatAn-uniri-iz-25-108 (Davor Dragi\v cevi\' c), and by the Japan Society for the Promotion of Science (JSPS) KAKENHI Grant number JP25K07092 (Masakazu Onitsuka).

\end{document}